\documentclass[11pt,a4paper]{amsart}
\usepackage{lscape}
\usepackage[latin1]{inputenc}
\usepackage{tikz}

\usetikzlibrary{shapes,arrows}

\tikzstyle{decision} = [diamond, draw, fill=white,
    text width=4.5em, text badly centered, node distance=3cm, inner sep=0pt]
\tikzstyle{block} = [rectangle, draw, fill=white,
    text width=12em, text centered, rounded corners, minimum height=4em]
\tikzstyle{line} = [draw, -latex']
\tikzstyle{cloud} = [draw, ellipse,fill=white, node distance=3cm,
    minimum height=2em]

 \theoremstyle{definition}

 \theoremstyle{remark}
 
 \theoremstyle{notation}
 
 \numberwithin{equation}{section}

\newcommand{\beq}{\begin{equation}}
\newcommand{\eeq}{\end{equation}}
\newcommand{\imn}{\mathbf{i}}

\begin{document}

\title[ The quadratic WDVV solution $E_8(a_1)$]
 {The quadratic WDVV solution $E_8(a_1)$ }

\author{Yassir Ibrahim Dinar }

\address{Faculty of Mathematical Sciences, University of Khartoum,  Sudan. Email: dinar@ictp.it.}
\email{dinar@ictp.it}


\begin{abstract}
We calculate explicitly  the quadratic solution to the WDVV equations corresponds to the quasi-Coxeter conjugacy class $E_8(a_1)$ using the associated classical $W$-algebra.
\end{abstract}
\maketitle

 A Frobenius algebra is a commutative associative algebra with unity $e$ and an invariant nondegenerate bilinear form $(.,.)$.
A \textbf{Frobenius manifold} is a manifold $M$ with
a smooth structure of Frobenius algebra on the tangent space $T_tM$ at
any point $t \in M $ with certain compatibility conditions \cite{DuRev}. Globally, we  require the metric $(.,.)$ to be flat and the unity vector field $e$ is constant with respect to it. In the flat  coordinates $(t^1,...,t^r)$ where $e={\partial\over \partial t^{r-1}}$ the compatibility conditions implies  that there exist a function $\mathbb{F}(t^1,...,t^r)$ such that
\[ \eta_{ij}=(\partial_{t^i},\partial_{t^j})=  \partial_{t^{r-1}}
\partial_{t^i}
\partial_{t^j} \mathbb{F}(t)\]
and the structure constants of the Frobenius algebra is given by
\[ C_{ij}^k=\eta^{kp}  \partial_{t^p}\partial_{t^i}\partial_{t^j} \mathbb{F}(t)\]
where $\eta^{ij}$ denote the inverse of the matrix $\eta_{ij}$.
In this work, we consider Frobenius manifolds where the quasihomogeneity condition takes the form
\begin{equation}
\sum_{i=1}^r d_i t^i \partial_{t^i} \mathbb{F}(t) = \left(3-d \right) \mathbb{F}(t);~~~~d_{r-1}=1.
\end{equation}
 This condition defines {\bf the degrees} $d_i$ and {\bf the charge} $d$ of  the Frobenius structure.
 If $\mathbb{F}(t)$ is an algebraic function we call $M$ an \textbf{algebraic Frobenius manifold}.
The  associativity of Frobenius
algebra implies the potential $\mathbb{F}(t)$ satisfies a system of partial differential equations known as \textbf{the  Witten-Dijkgraaf-Verlinde-Verlinde (WDVV) equations}:
\begin{equation} \label{frob}
 \partial_{t^i}
\partial_{t^j}
\partial_{t^k} \mathbb{F}(t)~ \eta^{kp} ~\partial_{t^p}
\partial_{t^q}
\partial_{t^n} \mathbb{F}(t) = \partial_{t^n}
\partial_{t^j}
\partial_{t^k} \mathbb{F}(t) ~\eta^{kp}~\partial_{t^p}
\partial_{t^q}
\partial_{t^i} \mathbb{F}(t).
  \end{equation}
In  topological field theory a solution to WDVV equations describes the a module space of two dimensional topological field theory \cite{wdvv}.

\textbf{Dubrovin conjecture} on classification of algebraic Frobenius manifolds and hence algebraic WDVV solutions, is stated as follows:
semisimple  irreducible algebraic Frobenius manifolds
with positive degrees $d_i$ correspond to quasi-Coxeter (primitive)  conjugacy
classes in irreducible Coxeter groups. A quasi-Coxeter conjugacy class in an irreducible Coxeter group is  a Conjugacy class  which has no representative in a proper Coxeter subgroup.

There are two major  results  support  the conjecture. First, the  conjecture arises from studying the algebraic solutions to associated equations of isomonodromic deformation of algebraic Frobenius manifolds \cite{DuRev},\cite{DMD}. It  leads to quasi-Coxeter conjugacy classes in Coxeter groups by considering the  classification of  finite orbits of the braid group action on tuple of reflections obtained in \cite{STF}.  Therefore, it remains the problem of constructing  all these algebraic  Frobenius manifolds. Second, Dubrovin constructed polynomial Frobenius structures on the orbit spaces of  Coxeter groups \cite{DCG}. Then  Hertling  \cite{HER} proved that  these are all possible \textbf{polynomial Frobenius manifolds}. The isomonodromic deformation of Polynomial Frobenius manifolds lead to   Coxeter conjugacy classes  \cite{DuRev}.

The classification of polynomial Frobenius manifolds reveals a relation between the order and eigenvalues of the  conjugacy class, and the charge and degrees of the corresponding Frobenius manifold. More precisely, If the order of the conjugacy class is $\kappa+1$ and the eigenvalues are $\exp {2\eta_i \pi \imn\over \kappa+1}$ then the charge of the Frobenius manifold is $\kappa-1\over \kappa+1$ and the degrees are $\eta_i+1\over \kappa+1$.  We depend on this \textbf{weak relation} in considering  a new examples of algebraic Frobenius manifolds.  In \cite{mypaper} we continue the work of  \cite{PAV} and we began to develop a construction of algebraic Frobenius manifolds using  classical $W$-algebras. This means we restrict ourself to conjugacy classes in Weyl groups. The examples  obtained  correspond, in the notations of \cite{CarClassif},  to the conjugacy classes $D_4(a_1)$ and $F_4(a_2)$.

In \cite{mypaper1} we uniform the construction of polynomial  Frobenius manifolds from classical $W$-algebras  associated to regular nilpotent orbits. In \cite{mypaper3} we extend this result and  we uniform the construction of algebraic Frobenius manifolds form the classical $W$-algebras associated to subregular nilpotent orbits in the Lie algebra of type $D_r$ where $r$ is even and $E_r$. These subregular nilpotent orbits correspond to  the regular quasi-Coxeter conjugacy classes $D_r(a_1)$ where $r$ is even and $E_r(a_1)$, respectively.

In this work we write explicitly the potential of the algebraic Frobenius manifold  $E_8(a_1)$. One of the main obstacles   was that  the minimal faithful matrix representation of the Lie algebra $E_8$ is the adjoint representation. This  means that the elements of the Lie algebra are represented by square matrices of dimension 248. We overcome this problem by constructing  instead   the Weyl-Chevalley normal form. Another problem arises was the following.  The Frobenius structure lives in a   hypersurface which  obtained by calculating  the restriction of the invariant polynomials of the adjoint action to Slodowy slice. We  avoid this calculation by returning  to the method  we used in \cite{mypaper} to obtain the algebraic Frobenius manifold $F_4(a_2)$.  Which  depends on the existence  of a  local Poisson  structure compatible with the classical $W$-algebras (bihamiltonian structure).

In this paper we will  not review  the rich and deep theory behind the construction of the WDVV solution  from classical $W$-algebras associated to nilpotent orbits. The interested reader may  consult the paper \cite{mypaper3} and \cite{mypaper} for details.

The resulting Frobenius manifold has the following potential $\mathbb F(t_1,t_2,...,t_8,Z)$ which  has 303 monomial. Here   $Z$ is a solution of a quadratic equation  where the  coefficients are polynomials in $t_1,...,t_6$ and $t_8$. We could not simply the potential $\mathbb{F}$ further. Any simplification comes with the price of losing the flat coordinates or having  very large quadratic equation for $Z$.  The quasihomogeneity reads
 \begin{equation}
 (\frac{1}{12} t_1 \partial_{t_1}+\frac{1}{3} t_2 \partial_{t_2}+\frac{1}{2} t_3 \partial_{t_3}+\frac{7}{12}
   t_4 \partial_{t_4}+\frac{3}{4} t_5 \partial_{t_5}+\frac{5}{6} t_6 \partial_{t_6}+t_7 \partial_{t_7}+\frac{1}{4}
   t_8 \partial_{t_8})\mathbb F={25\over 12} \mathbb F.
   \end{equation}
 In the end we would like to add that we are very  impressed  by writing this potential in a paper. But we believe that, this potential   gives an idea about how complicated algebraic WDVV solutions could be. We write this potential  in   a Mathematica notebook file and we will  uploaded it in the  sub-directory  of the arXiv submission. It can also be found on the  homepage http://staffcv.uofk.edu/FMS/Dept-of-Pure-Math/Yassir-Ibrahim/.

\begin{landscape}
\footnotesize{
\begin{eqnarray*}
F&=&\frac{275747251366536586569731516102824113952687921718886400000000
   t_1^{25}}{262919712977108888812209859761}-\frac{329401897854852815100761902648244065166950400000000
   \sqrt{\frac{2}{134589}} t_8
   t_1^{22}}{2736152181450029070729}
   \\ & &-\frac{1945051421153779792882712347780124853351219200000
   \sqrt{\frac{288230}{3289}} t_2
   t_1^{21}}{970315921277476588361547}+\frac{60257490748451656689759529090416640000000000
   \sqrt{\frac{5510}{69069}} Z
   t_1^{20}}{51511170636379284831}
      \\ & &+\frac{51276965550447936573186450702116293967872000000 t_8^2
   t_1^{19}}{8491244341173937384041}+\frac{8609664200367436587684526517904002252800000
   \sqrt{\frac{21793}{221}} t_3
   t_1^{19}}{10845082617828469429419}
      \\ & &+\frac{237507977941170664487849007390455234560000
   \sqrt{\frac{19}{19499}} t_4
   t_1^{18}}{380778653527744909437}-\frac{1158908336560292796842981341375391334400000
   \sqrt{\frac{144115}{2619309}} t_2 t_8
   t_1^{18}}{51327530099351194689}
      \\ & &+\frac{734109387995764004804539557492788341768192000 t_2^2
   t_1^{17}}{6945927852281967485937}-\frac{274536606777826836500216545280000000000
   \sqrt{\frac{95}{4301}} t_8 Z
   t_1^{17}}{57221326755129537}
      \\ & &-\frac{5452777995210797138612019621629054156800000
   \sqrt{\frac{2}{134589}} t_8^3
   t_1^{16}}{2384044104068998371}+\frac{2055990113756671264610881296662528000 \sqrt{\frac{190}{383801}}
   t_5 t_1^{16}}{305692228025763579}
      \\ & &+\frac{87138535415522897242687767838720000 \sqrt{\frac{1263994}{21}}
   t_3 t_8 t_1^{16}}{35723692959645537}+\frac{10638293512640789914383391129600000000
   \sqrt{\frac{43993}{21}} t_2 Z
   t_1^{16}}{14380803006794963739}
      \\ & &-\frac{26661242077408439262555938492514304000
   \sqrt{\frac{288230}{3289}} t_2 t_8^2
   t_1^{15}}{240252506611604487}+\frac{237239573650736748213584251337886924800 \sqrt{\frac{12710}{4301}}
   t_2 t_3 t_1^{15}}{11967241662997207587}
   \\ & &+\frac{552547343072105402364174348478054400
   \sqrt{\frac{58}{357}} t_6 t_1^{15}}{6268555156808061117}+\frac{19572106164686166880069029889835008000
   \sqrt{\frac{38}{9080799}} t_4 t_8
   t_1^{15}}{26339975807916771}
   \\ & &+\frac{649005721094941689805827190292480 \sqrt{\frac{12710}{55913}} t_2
   t_4 t_1^{14}}{3591268931403153}+\frac{9634301289821016796784088857742147584000
   \sqrt{\frac{2}{134589}} t_2^2 t_8
   t_1^{14}}{367423962740453619}
   \\ & &+\frac{1131978190761833389673676800000000 \sqrt{\frac{4370}{87087}}
   t_8^2 Z t_1^{14}}{91039221906633}-\frac{172298337875627126439280640000000
   \sqrt{\frac{332630}{90321}} t_3 Z
   t_1^{14}}{635305073153751}
   \\ & &+\frac{25908989728497747627499196136816640000 t_8^4
   t_1^{13}}{2309947939218369051}-\frac{27498013460678823586752388578267889664
   \sqrt{\frac{288230}{3289}} t_2^3
   t_1^{13}}{22993904272896760203}
   \\ & &+\frac{7756364878391254176136276023758028800 t_3^2
   t_1^{13}}{1407624273750994299}-\frac{69171011505030796144541696000 \sqrt{\frac{283309}{17}} t_3 t_8^2
   t_1^{13}}{1889703018819}
   \\ & &+\frac{37858410608405201801558622208000 \sqrt{\frac{2185}{2245886643}} t_5
   t_8 t_1^{13}}{3178802752587}-\frac{4753057596569024177635328000000 \sqrt{\frac{10}{39056516499}} t_4
   Z t_1^{13}}{59167157049}
   \\ & &+\frac{12117169859674321505484800000000 \sqrt{\frac{3034}{221}} t_2 t_8 Z
   t_1^{13}}{175548954964383}+\frac{682332614278050128408269448806400 \sqrt{\frac{3314645}{113883}} t_2
   t_8^3 t_1^{12}}{2630725180621083}
    \end{eqnarray*}}
\footnotesize{\begin{eqnarray*}
& &-\frac{581596580070975576125721070796800 \sqrt{\frac{19}{19499}} t_4
   t_8^2 t_1^{12}}{2168478895654071}-\frac{59277167362036445126003261440 t_3 t_4
   t_1^{12}}{70589452572639 \sqrt{13}}+\frac{56181243169575977303136010240 t_2 t_5
   t_1^{12}}{23529817524213 \sqrt{13}}
   \\ & &-\frac{245825012333495379886125130711040
   \sqrt{\frac{184295}{69069}} t_2 t_3 t_8
   t_1^{12}}{275940587329407}+\frac{236615066302532511259741388800 t_6 t_8 t_1^{12}}{23529817524213
   \sqrt{13}}+\frac{3748908861782265246305484800000 \sqrt{\frac{5510}{69069}} t_2^2 Z
   t_1^{12}}{93898278236763}
   \\ & &+\frac{3775956970890476910732497649664 t_4^2
   t_1^{11}}{906846213331879821}+\frac{5002923413852734025203121190338560 t_2^2 t_8^2
   t_1^{11}}{28745698069975887}+\frac{35514283134452206559049209085952 \sqrt{\frac{21793}{221}} t_2^2
   t_3 t_1^{11}}{2196577671230943}
   \\ & &+\frac{24634735902975255837197467648 \sqrt{\frac{835867}{5870865}} t_2
   t_6 t_1^{11}}{59121138149073}-\frac{48882242654346106005782265856 \sqrt{\frac{6355}{154077}} t_2 t_4
   t_8 t_1^{11}}{9688472554761}-\frac{62651926800793631457280000000 \sqrt{\frac{2185}{187}} t_8^3 Z
   t_1^{11}}{67454491810797}
   \\ & &-\frac{411448891669756248064000000 \sqrt{\frac{29}{414141}} t_5 Z
   t_1^{11}}{348506614071}+\frac{78500117489624547328000000 \sqrt{\frac{5735}{3289}} t_3 t_8 Z
   t_1^{11}}{54286954029}-\frac{30579366843514835600438422667264000 \sqrt{\frac{2}{134589}} t_8^5
   t_1^{10}}{10304778647825487}
   \\ & &+\frac{22232868668243672938263347200 \sqrt{\frac{43586}{609}} t_3 t_8^3
   t_1^{10}}{17236382080743}-\frac{125348126861275384690769920 \sqrt{\frac{161690}{451}} t_5 t_8^2
   t_1^{10}}{40040071035111}-\frac{189760322042691934295687168 \sqrt{370481} t_2^2 t_4
   t_1^{10}}{936498700932795}
   \\ & &-\frac{14880867845386572962725888 \sqrt{\frac{130}{5466571}} t_3 t_5
   t_1^{10}}{40539870063}-\frac{62830748616325803253808902438912 \sqrt{\frac{2449955}{154077}} t_2^3 t_8
   t_1^{10}}{20677539028683537}-\frac{95011219805712869211813969920 \sqrt{\frac{754}{357}} t_3^2 t_8
   t_1^{10}}{40093931492307}
   \\ & &-\frac{161203054336679084032000000 \sqrt{\frac{1517}{609}} t_2 t_8^2 Z
   t_1^{10}}{1105051246299}+\frac{29618300908091893350400000 \sqrt{\frac{700321}{4641}} t_2 t_3 Z
   t_1^{10}}{13643339572863}-\frac{33755050520538555351040000 \sqrt{\frac{95}{55913}} t_6 Z
   t_1^{10}}{593188677081}
   \\ & &+\frac{2165520482472401305600000 \sqrt{\frac{115}{12617}} t_4 t_8 Z
   t_1^{10}}{106466731371}+\frac{379599576880637346221842511367847936 t_2^4
   t_1^9}{248022271919052355725}-\frac{31529410268465987578882949120 \sqrt{\frac{72922190}{13}} t_2
   t_8^4 t_1^9}{28038583762687953}
   \\ & &+\frac{197407620359075269146711162880 \sqrt{\frac{38}{9080799}} t_4
   t_8^3 t_1^9}{31050480357351}-\frac{3160121833455978507411980288 \sqrt{\frac{57646}{16445}} t_2 t_3^2
   t_1^9}{7148530421109}+\frac{217861791934905227836129280 \sqrt{\frac{12710}{4301}} t_2 t_3 t_8^2
   t_1^9}{306529366689}
   \\ & &-\frac{1260293666811736204771328 \sqrt{\frac{2}{10353}} t_6 t_8^2
   t_1^9}{1244674431}-\frac{2052533495915389374169088 \sqrt{\frac{46}{1188385}} t_4 t_5
   t_1^9}{2888484111513}-\frac{1052429798275695127298048 \sqrt{\frac{1263994}{273}} t_3 t_6
   t_1^9}{168501074211855}
   \\ & &-\frac{109423644938076252274688 \sqrt{\frac{2}{10353}} t_3 t_4 t_8
   t_1^9}{86837751}-\frac{185359867371260393553920 \sqrt{\frac{2}{10353}} t_2 t_5 t_8
   t_1^9}{339456663}-\frac{712342263971184640000 \sqrt{\frac{410533}{609}} t_2 t_4 Z
   t_1^9}{280809205677}
   \\ & &+\frac{34160490430294733619200000 \sqrt{\frac{1615}{253}} t_2^2 t_8 Z
   t_1^9}{23051882975121}-\frac{54459877768149145360567107584 \sqrt{\frac{34}{7917}} t_2^2 t_8^3
   t_1^8}{13919180849847}+\frac{169398250083998769348608 \sqrt{\frac{993922}{715}} t_2 t_4 t_8^2
   t_1^8}{3129118175547}
   \\ & &+\frac{4431201496798079746048 \sqrt{\frac{57646}{1265}} t_2 t_3 t_4
   t_1^8}{22909841955}-\frac{5205387859279057518592 \sqrt{\frac{57646}{1265}} t_2^2 t_5
   t_1^8}{129822437745}-\frac{145162730796647603765248 \sqrt{\frac{38}{698523}} t_4 t_6
   t_1^8}{1020033439425}
   \\ & &-\frac{11548884065344094489018368 \sqrt{\frac{34}{7917}} t_4^2 t_8
   t_1^8}{218605116394665}+\frac{105526751443021144048369664 \sqrt{\frac{1263994}{21}} t_2^2 t_3 t_8
   t_1^8}{128394390929805}-\frac{1890182858062194802688 \sqrt{\frac{57646}{1265}} t_2 t_6 t_8
   t_1^8}{7081223877}
    \end{eqnarray*}}
\footnotesize{\begin{eqnarray*}
    & &+\frac{2391235564049268736000000 \sqrt{\frac{48070}{7917}} t_8^4 Z
   t_1^8}{1180522241703}+\frac{22503223070956655771648000 \sqrt{\frac{43993}{21}} t_2^3 Z
   t_1^8}{1931124907440729}+\frac{7435342160099246080000 \sqrt{\frac{71630}{5313}} t_3^2 Z
   t_1^8}{133202430207}
   \\ & &-\frac{28365429486387200000 \sqrt{\frac{7650490}{3927}} t_3 t_8^2 Z
   t_1^8}{374272731}+\frac{10778863204827136000000 \sqrt{\frac{2}{25789}} t_5 t_8 Z
   t_1^8}{8904209391}+\frac{21434549150503361944170895769600 t_8^6
   t_1^7}{4283352991279460763}
   \\ & &-\frac{4287547480979049202319360 \sqrt{\frac{283309}{17}} t_3 t_8^4
   t_1^7}{46827382254219}+\frac{469684633100594953600000 \sqrt{\frac{283309}{17}} t_3^3
   t_1^7}{21445591263327}+\frac{74052171417443891150848 \sqrt{\frac{11560835}{424473}} t_5 t_8^3
   t_1^7}{4822753654089}
   \\ & &-\frac{265329980759757291520 \sqrt{\frac{7790}{121693}} t_2 t_4^2
   t_1^7}{103715653041}+\frac{2828941945374762065199104 t_5^2
   t_1^7}{10944541213524765}-\frac{15887825476612819786628661248 \sqrt{\frac{57646}{16445}} t_2^3 t_8^2
   t_1^7}{1451802932827065}
   \\ & &+\frac{5706043389286903615127552 t_3^2 t_8^2
   t_1^7}{757153734813}-\frac{381865658996990733088897719296 \sqrt{\frac{2542}{21505}} t_2^3 t_3
   t_1^7}{7089794240604075}-\frac{328941634965862444714172416 \sqrt{\frac{58}{357}} t_2^2 t_6
   t_1^7}{2063167636154625}
   \\ & &-\frac{1287399644481804608602112 \sqrt{\frac{43586}{7917}} t_2^2 t_4 t_8
   t_1^7}{88426535222025}-\frac{1875868881215533088768 \sqrt{\frac{667}{1468005}} t_3 t_5 t_8
   t_1^7}{838486971}+\frac{18234492407722803200000 \sqrt{\frac{51578}{13}} t_2 t_8^3 Z
   t_1^7}{1070706219219}
   \\ & &-\frac{90769374356439040000 \sqrt{\frac{43010}{9080799}} t_4 t_8^2 Z
   t_1^7}{1185721173}+\frac{12332795428864000 \sqrt{\frac{5510}{5313}} t_3 t_4 Z
   t_1^7}{12405393}-\frac{6166397714432000 \sqrt{\frac{5510}{5313}} t_2 t_5 Z
   t_1^7}{12405393}-\frac{256391108968939520000 \sqrt{48298} t_2 t_3 t_8 Z
   t_1^7}{257647607217}
   \\ & &+\frac{30831988572160000 \sqrt{\frac{5510}{5313}} t_6 t_8 Z
   t_1^7}{12405393}+\frac{46404673915505485286998016 \sqrt{\frac{32623085}{11571}} t_2 t_8^5
   t_1^6}{2126382895583037}-\frac{64219603583178549716058112 \sqrt{\frac{323}{1147}} t_4 t_8^4
   t_1^6}{4302216556179633}
   \\ & &-\frac{1871220704220541288448 \sqrt{\frac{321563}{39585}} t_2 t_3 t_8^3
   t_1^6}{2299103919}   +\frac{3722731443539978223616 t_6 t_8^3 t_1^6}{54269545995
   \sqrt{13}}+\frac{698055308833678426112 t_3 t_4 t_8^2 t_1^6}{18089848665
   \sqrt{13}}+\frac{106374023157810987008 t_2 t_5 t_8^2 t_1^6}{2584264095
   \sqrt{13}}
   \\ & &+\frac{296437133234144764850176 \sqrt{\frac{216070}{3289}} t_2^3 t_4
   t_1^6}{221354238402297}   +\frac{183825589760000 \sqrt{\frac{21793}{17}} t_3^2 t_4
   t_1^6}{51212007}-\frac{16767557595947008 \sqrt{\frac{21793}{17}} t_2 t_3 t_5
   t_1^6}{55665225}-\frac{12566025865360769024 \sqrt{\frac{551}{685084785}} t_5 t_6
   t_1^6}{1396932075}
   \\ & &+\frac{4508542535356502346462936064 \sqrt{\frac{34}{7917}} t_2^4 t_8
   t_1^6}{686267529090507}   -\frac{3056148141650068087808 \sqrt{\frac{835867}{451605}} t_2 t_3^2 t_8
   t_1^6}{5672117745}-\frac{6213465218395906506752 \sqrt{\frac{23}{553437885}} t_4 t_5 t_8
   t_1^6}{22977914085}
   \\ & &-\frac{6882616087455776768 \sqrt{\frac{21793}{17}} t_3 t_6 t_8
   t_1^6}{21765102975}+\frac{83863008916275200 \sqrt{\frac{4370}{87087}} t_4^2 Z
   t_1^6}{15026127831}
   +\frac{16457073999942434816000 \sqrt{\frac{190}{2003001}} t_2^2 t_8^2 Z
   t_1^6}{1594590543} -\frac{42878039407831859200 \sqrt{\frac{332630}{90321}} t_2^2 t_3 Z
   t_1^6}{15055038999}
    \\ & &-\frac{1489840227792486400 \sqrt{\frac{3034}{17}} t_2 t_6 Z
   t_1^6}{165605428989}-\frac{104828761145344000 \sqrt{\frac{48298}{13}} t_2 t_4 t_8 Z
   t_1^6}{11829636819}   -\frac{1540964400458670051556291472 \sqrt{\frac{57646}{16445}} t_2^5
   t_1^5}{52491486120434361}
   \\ & &+\frac{2950962167502149088348848128 t_2^2 t_8^4
   t_1^5}{26173583482176585}-\frac{397577030236771188736 \sqrt{\frac{321563}{3045}} t_2 t_4 t_8^3
   t_1^5}{3559578876945}+\frac{23921384687434934326958624 t_2^2 t_3^2
   t_1^5}{172706092106625}+\frac{96349964288000 \sqrt{\frac{323}{14911}} t_3 t_4^2
   t_1^5}{2351349}
   \\ & &+\frac{1607790932517785666048 t_6^2
   t_1^5}{255662227445625}+\frac{32316116749389439434752 t_4^2 t_8^2
   t_1^5}{5900445457665525}-\frac{34013318595020812490752 \sqrt{\frac{1147}{4199}} t_2^2 t_3 t_8^2
   t_1^5}{164569921425}
       \end{eqnarray*}}
\footnotesize{\begin{eqnarray*}
   & &+\frac{235811979254300672 \sqrt{\frac{28823}{170255085}} t_2 t_6 t_8^2
   t_1^5}{1131975}+\frac{1953069824573440 \sqrt{\frac{19}{253487}} t_2 t_4 t_5
   t_1^5}{2351349}-\frac{30272806725074369536 \sqrt{\frac{36859}{26565}} t_2 t_3 t_6
   t_1^5}{39062219625}+\frac{785230641485037568 \sqrt{\frac{28823}{170255085}} t_2 t_3 t_4 t_8
   t_1^5}{7923825}
   \\ & &+\frac{2898792624381952 \sqrt{\frac{144115}{34051017}} t_2^2 t_5 t_8
   t_1^5}{316953}-\frac{250359054902755328 \sqrt{\frac{19}{253487}} t_4 t_6 t_8
   t_1^5}{293918625}-\frac{181117057787494400000 \sqrt{408595} t_8^5 Z
   t_1^5}{54522638348283}+\frac{6616808161280000 \sqrt{18862415} t_3 t_8^3 Z
   t_1^5}{5364575811}
   \\ & &-\frac{587345111941120000 \sqrt{\frac{13}{923853}} t_5 t_8^2 Z
   t_1^5}{61388649}-\frac{8765615487877120 \sqrt{\frac{194990}{2003001}} t_2^2 t_4 Z
   t_1^5}{135069363}+\frac{402481291264000 \sqrt{\frac{551}{453747}} t_3 t_5 Z
   t_1^5}{1802493}+\frac{205051967041083392000 \sqrt{\frac{51578}{13}} t_2^3 t_8 Z
   t_1^5}{2805255600147}
   \\ & &-\frac{352308820824064000 \sqrt{\frac{95}{4301}} t_3^2 t_8 Z
   t_1^5}{49322763}-\frac{2158936234052031567626240000 \sqrt{\frac{442}{609}} t_8^7
   t_1^4}{2645124258409160577}+\frac{8695021482453436416000 \sqrt{\frac{6882}{3857}} t_3 t_8^5
   t_1^4}{3551282957551}-\frac{336957034612426342400 \sqrt{\frac{48070}{1517}} t_5 t_8^4
   t_1^4}{74246590205543}
   \\ & &+\frac{26965819478835200 \sqrt{\frac{323}{1147}} t_4^3
   t_1^4}{54760218525441}+\frac{932722984197401182208 \sqrt{\frac{26949505}{14007}} t_2^3 t_8^3
   t_1^4}{8616635764191}-\frac{25629396450036121600 \sqrt{\frac{754}{357}} t_3^2 t_8^3
   t_1^4}{10133087643}-\frac{5833918962368512 \sqrt{370481} t_2^2 t_4 t_8^2
   t_1^4}{45048212937}
   \\ & &+\frac{10976384438456320 \sqrt{\frac{32890}{21607}} t_3 t_5 t_8^2
   t_1^4}{235653201}-\frac{13241958273039104 t_2^2 t_3 t_4 t_1^4}{21724443
   \sqrt{13}}-\frac{34145463252284672 t_2^3 t_5 t_1^4}{386026641 \sqrt{13}}+\frac{31195609952000
   \sqrt{\frac{7030}{10373}} t_3^2 t_5 t_1^4}{928557}-\frac{44167651011584 \sqrt{\frac{8897}{1430715}}
   t_2 t_4 t_6 t_1^4}{3463317}
   \\ & &-\frac{37664594331520000 \sqrt{\frac{1263994}{21}} t_3^3 t_8
   t_1^4}{19285200333}-\frac{16549999230976 \sqrt{\frac{1522945}{247863}} t_2 t_4^2 t_8
   t_1^4}{238350333}+\frac{442283885249822720 \sqrt{\frac{26}{10353}} t_5^2 t_8
   t_1^4}{380447348379}-\frac{2792998290203006336 \sqrt{\frac{405449}{31395}} t_2^3 t_3 t_8
   t_1^4}{726561927}
   \\ & &+\frac{144777055201009664 t_2^2 t_6 t_8 t_1^4}{363648285
   \sqrt{13}}+\frac{269622764011015604416 \sqrt{\frac{5510}{69069}} t_2^4 Z
   t_1^4}{928872009351}-\frac{61124597799485440000 \sqrt{\frac{1517}{609}} t_2 t_8^4 Z
   t_1^4}{243599827653}+\frac{493179660795904000 \sqrt{\frac{1265}{1147}} t_4 t_8^3 Z
   t_1^4}{164288255859}
   \\ & &-\frac{2625378232102400 \sqrt{\frac{43993}{21}} t_2 t_3^2 Z
   t_1^4}{2860556391}+\frac{23320253440000 \sqrt{\frac{700321}{4641}} t_2 t_3 t_8^2 Z
   t_1^4}{349461}-\frac{5618892800 \sqrt{\frac{28405}{187}} t_6 t_8^2 Z
   t_1^4}{8127}+\frac{19643649228800 \sqrt{\frac{19}{171062619}} t_4 t_5 Z
   t_1^4}{340659}
   \\ & &+\frac{48259054827520 \sqrt{\frac{5735}{253}} t_3 t_6 Z
   t_1^4}{124372017}-\frac{11237785600 \sqrt{\frac{2185}{2431}} t_3 t_4 t_8 Z
   t_1^4}{3483}-\frac{112377856000 \sqrt{\frac{2185}{2431}} t_2 t_5 t_8 Z
   t_1^4}{24381}-\frac{1297036772272472075468800 \sqrt{\frac{3838010}{247}} t_2 t_8^6
   t_1^3}{1291804870385869119}
   \\ & &+\frac{109310473017890278604800 \sqrt{\frac{6}{57511727}} t_4 t_8^5
   t_1^3}{25097722414289}+\frac{1008644456653419520 \sqrt{\frac{3215630}{17}} t_2 t_3 t_8^4
   t_1^3}{830737873791}-\frac{78331776737280 \sqrt{\frac{6}{3451}} t_6 t_8^4
   t_1^3}{556549}-\frac{323750049601694950 \sqrt{\frac{12710}{4301}} t_2 t_3^3
   t_1^3}{1820354067}
   \\ & &-\frac{4905667695554560 \sqrt{\frac{2}{10353}} t_3 t_4 t_8^3
   t_1^3}{16139921}-\frac{11919985887293440 \sqrt{\frac{2}{10353}} t_2 t_5 t_8^3
   t_1^3}{48419763}+\frac{33332406263111680 t_2^2 t_4^2 t_1^3}{492237781893}-\frac{11695614464
   \sqrt{\frac{56810}{16687}} t_2 t_5^2 t_1^3}{564417}+\frac{324321642022604931470608 t_2^4 t_8^2
   t_1^3}{1378060561635849}
   \\ & &-\frac{4735581833401184 \sqrt{\frac{2169310}{437}} t_2 t_3^2 t_8^2
   t_1^3}{706959603}+\frac{457125753110528 \sqrt{\frac{2530}{21607}} t_4 t_5 t_8^2
   t_1^3}{5350785969}+\frac{10308305792 \sqrt{\frac{16431922}{21}} t_3 t_6 t_8^2
   t_1^3}{782901}+\frac{3337040790377704912907 \sqrt{\frac{370481}{13}} t_2^4 t_3
   t_1^3}{1035335031961230}
   \\ & &-\frac{702361600 \sqrt{\frac{56810}{16687}} t_3 t_4 t_5
   t_1^3}{33201}+\frac{35772423689400752 \sqrt{\frac{99468173}{49335}} t_2^3 t_6
   t_1^3}{1071775395405}-\frac{4808120870133100 \sqrt{\frac{58}{357}} t_3^2 t_6
   t_1^3}{317839599}+\frac{10524737455707266816 \sqrt{\frac{1271}{770385}} t_2^3 t_4 t_8
   t_1^3}{199055403963}
   \\ & &-\frac{11552000 \sqrt{\frac{16431922}{21}} t_3^2 t_4 t_8
   t_1^3}{2709}+\frac{33162496 \sqrt{\frac{16431922}{21}} t_2 t_3 t_5 t_8
   t_1^3}{111843}-\frac{92015148032 \sqrt{\frac{11362}{83435}} t_5 t_6 t_8
   t_1^3}{564417}-\frac{2989531521228800 \sqrt{408595} t_2^2 t_8^3 Z
   t_1^3}{73646459523}
          \end{eqnarray*}}
\footnotesize{\begin{eqnarray*}
    & &+\frac{231261798400 \sqrt{\frac{410533}{609}} t_2 t_4 t_8^2 Z
   t_1^3}{18820971}-\frac{2177320960 \sqrt{\frac{43993}{273}} t_2 t_3 t_4 Z
   t_1^3}{266409}-\frac{1088660480 \sqrt{\frac{43993}{273}} t_2^2 t_5 Z
   t_1^3}{88803}+\frac{45906354176 \sqrt{\frac{5}{3772483}} t_4 t_6 Z
   t_1^3}{11583}
     \\ & &-\frac{3515179335680 \sqrt{\frac{37145}{11}} t_4^2 t_8 Z
   t_1^3}{189268499043}+\frac{742348327697920 \sqrt{\frac{5735}{3289}} t_2^2 t_3 t_8 Z
   t_1^3}{21869757}-\frac{21773209600 \sqrt{\frac{43993}{273}} t_2 t_6 t_8 Z
   t_1^3}{266409}-\frac{155540097158692339850240 \sqrt{\frac{102}{2639}} t_2^2 t_8^5
   t_1^2}{4061503837595113}
     \\ & &+\frac{11540725369823232 \sqrt{\frac{54665710}{13}} t_2 t_4 t_8^4
   t_1^2}{144165986891381}-\frac{405070505245294592 \sqrt{\frac{34}{7917}} t_4^2 t_8^3
   t_1^2}{28913863781153}+\frac{90733129640594848 \sqrt{\frac{2294}{11571}} t_2^2 t_3 t_8^3
   t_1^2}{1237208427}-\frac{360361606016 \sqrt{\frac{767602}{95}} t_2 t_6 t_8^3
   t_1^2}{43915599}
     \\ & &-\frac{10159893376 \sqrt{\frac{767602}{95}} t_2 t_3 t_4 t_8^2
   t_1^2}{6273657}-\frac{710272 \sqrt{72922190} t_2^2 t_5 t_8^2 t_1^2}{48633}+\frac{32689198592
   \sqrt{\frac{266}{99789}} t_4 t_6 t_8^2 t_1^2}{243165}-\frac{164315483619570550 \sqrt{370481} t_2^4
   t_4 t_1^2}{350184747525327}+\frac{240158860}{161} \sqrt{\frac{12710}{55913}} t_2 t_3^2 t_4
   t_1^2
   \\ & &+\frac{54924677120 \sqrt{\frac{4370}{16687}} t_4^2 t_5 t_1^2}{1166234043}+\frac{35601655004
   \sqrt{\frac{33046}{21505}} t_2^2 t_3 t_5 t_1^2}{73899} -\frac{47189920 \sqrt{\frac{58}{4641}} t_3 t_4
   t_6 t_1^2}{2277}+\frac{4668324928 \sqrt{\frac{58}{4641}} t_2 t_5 t_6
   t_1^2}{193545}   -\frac{153751727659501843325 \sqrt{\frac{2449955}{154077}} t_2^5 t_8
   t_1^2}{94407173314038}
   \\ & & +\frac{928913367501436 \sqrt{\frac{406}{663}} t_2^2 t_3^2 t_8
   t_1^2}{58346055}-\frac{13129464512 \sqrt{\frac{406}{663}} t_6^2 t_8 t_1^2}{967725}+\frac{35542629376
   \sqrt{\frac{38}{698523}} t_2 t_4 t_5 t_8 t_1^2}{48633}
   +\frac{3915730205848 \sqrt{\frac{2542}{279565}}
   t_2 t_3 t_6 t_8 t_1^2}{369495}
   \\ & &+\frac{5705857874083840000 \sqrt{\frac{32890}{11571}} t_8^6 Z
   t_1^2}{37213864269463} -\frac{1704175616000 \sqrt{\frac{2901910}{10353}} t_3 t_8^4 Z
   t_1^2}{23931607}+\frac{1489449488384000 \sqrt{\frac{2}{25789}} t_5 t_8^3 Z
   t_1^2}{1501012653}
   +\frac{1105295360 \sqrt{\frac{1647490}{231}} t_5^2 Z
   t_1^2}{1263944679}
   \\ & &-\frac{317637724130560 \sqrt{\frac{1517}{609}} t_2^3 t_8^2 Z
   t_1^2}{840886371}  +\frac{5534963200 \sqrt{\frac{18122390}{21}} t_3^2 t_8^2 Z
   t_1^2}{3740527}  +\frac{7370788883912 \sqrt{\frac{11905457}{273}} t_2^3 t_3 Z
   t_1^2}{1832094693}-\frac{21000643480928 \sqrt{\frac{323}{16445}} t_2^2 t_6 Z
   t_1^2}{79656291}\\ & &+\frac{135784610816 \sqrt{\frac{5735}{253}} t_2^2 t_4 t_8 Z
   t_1^2}{73970793}-\frac{8435148800 \sqrt{\frac{494}{1271}} t_3 t_5 t_8 Z
   t_1^2}{260967}+\frac{1453995723110575582080000000 t_8^8
   t_1}{146232167793805896952009}+\frac{4916728097717507366510447 t_2^6
   t_1}{21741928103909782080}
     \\ & &-\frac{3098936846720000000 \sqrt{\frac{253487}{19}} t_3 t_8^6
   t_1}{4949443947838579}+\frac{40567900538880000 \sqrt{\frac{838695}{5851069}} t_5 t_8^5
   t_1}{26828599822171}+\frac{255540004042051625 t_3^4 t_1}{534432838464}-\frac{2248830784740736064
   \sqrt{\frac{3838010}{247}} t_2^3 t_8^4 t_1}{1744905565478439}
     \\ & &+\frac{11166926928300000 t_3^2 t_8^4
   t_1}{54109363427}+\frac{1154280198676352 \sqrt{\frac{2294}{150423}} t_2^2 t_4 t_8^3
   t_1}{33900751209}-\frac{314691520000 \sqrt{\frac{3795}{258013}} t_3 t_5 t_8^3
   t_1}{9461333}-\frac{370421080518379193 \sqrt{\frac{1517}{624910}} t_2^3 t_3^2
   t_1}{7623051345}
     \\ & &-\frac{392768 \sqrt{\frac{216070}{253}} t_2 t_3 t_4^2 t_1}{44109}+\frac{4283474000
   \sqrt{\frac{247}{19499}} t_3 t_5^2 t_1}{2571233}-\frac{468277715963 \sqrt{\frac{28823}{32890}} t_2
   t_6^2 t_1}{868049325}+\frac{1}{2} t_7^2 t_1+\frac{1120775995625 \sqrt{\frac{14911}{323}} t_3^3 t_8^2
   t_1}{78551067}-\frac{105422908928 \sqrt{\frac{394174}{2405}} t_2 t_4^2 t_8^2
   t_1}{22863297327}
     \\ & &+\frac{412560582720000 t_5^2 t_8^2 t_1}{9302310655463}+\frac{135254310447593
   \sqrt{\frac{475354}{115}} t_2^3 t_3 t_8^2 t_1}{25808743233}-\frac{1043811292 \sqrt{\frac{238}{87}}
   t_2^2 t_6 t_8^2 t_1}{58305}-\frac{490960 \sqrt{\frac{216070}{253}} t_2^2 t_4 t_5
   t_1}{44109}-\frac{320891456 \sqrt{\frac{238}{87}} t_4^2 t_6 t_1}{185507751}
     \\ & &-\frac{19489949684819   \sqrt{\frac{631997}{546}} t_2^2 t_3 t_6 t_1}{79860537900}-\frac{19214964674560
   \sqrt{\frac{266}{1297257}} t_4^3 t_8 t_1}{609865163583}-\frac{82122873884 \sqrt{\frac{34}{609}} t_2^2
   t_3 t_4 t_8 t_1}{578565}-\frac{47088135100 \sqrt{\frac{34}{609}} t_2^3 t_5 t_8
   t_1}{1504269}+\frac{4901000 \sqrt{\frac{335264215}{14637}} t_3^2 t_5 t_8 t_1}{86989}
     \\ & &+\frac{18834016   \sqrt{\frac{43214}{1265}} t_2 t_4 t_6 t_8 t_1}{220545}+\frac{8739209194496000 \sqrt{\frac{51578}{13}}
   t_2 t_8^5 Z t_1}{16443335374879}-\frac{20502543564800 \sqrt{\frac{43010}{9080799}} t_4 t_8^4 Z
   t_1}{2875211641}-\frac{53859449600 \sqrt{\frac{2542}{19}} t_2 t_3 t_8^3 Z
   t_1}{5008941}+\frac{1295360 \sqrt{\frac{48070}{609}} t_6 t_8^3 Z t_1}{5547}
\end{eqnarray*}}
\footnotesize{\begin{eqnarray*}
 & &+\frac{235520   \sqrt{\frac{48070}{609}} t_3 t_4 t_8^2 Z t_1}{1849}+\frac{117760 \sqrt{\frac{48070}{609}} t_2 t_5
   t_8^2 Z t_1}{1849}-\frac{17632}{153} \sqrt{\frac{4324190}{90321}} t_2 t_3 t_5 Z
   t_1+\frac{57621376 \sqrt{\frac{26}{25789}} t_5 t_6 Z t_1}{52173}+\frac{2392274960
   \sqrt{\frac{39442}{17}} t_2 t_3^2 t_8 Z t_1}{1059219}
     \\ & &-\frac{205844480 \sqrt{\frac{38}{1271}} t_4
   t_5 t_8 Z t_1}{4290927}-\frac{35264}{51} \sqrt{\frac{4324190}{90321}} t_3 t_6 t_8 Z
   t_1+\frac{93968346863818240000 \sqrt{\frac{97869255}{3857}} t_2
   t_8^7}{16752453636591350321}-\frac{380282459448883200000 \sqrt{\frac{17}{21793}} t_4
   t_8^6}{594640337161749277}
   \\ & &-\frac{109833811208000 \sqrt{\frac{62704785}{203}} t_2 t_3
   t_8^5}{10773219728369} +\frac{1441585094400 \sqrt{13} t_6 t_8^5}{54629021179}+\frac{9610567296000
   \sqrt{13} t_3 t_4 t_8^4}{382403148253}+\frac{2402641824000 \sqrt{13} t_2 t_5
   t_8^4}{382403148253}+\frac{794569664 \sqrt{\frac{32062}{22165}} t_2
   t_4^3}{17843016789}
   \\ & &+\frac{574284891088122335 \sqrt{\frac{34}{7917}} t_2^4
   t_8^3}{59769362364309}+\frac{1398279025 \sqrt{\frac{55651145}{6783}} t_2 t_3^2
   t_8^3}{10574431}
    -\frac{174737587200 \sqrt{\frac{49335}{258013}} t_4 t_5
   t_8^3}{76694598889}-\frac{3946800 \sqrt{\frac{1147}{323}} t_3 t_6 t_8^3}{12943}-\frac{125436480
   \sqrt{\frac{19}{19499}} t_4 t_5^2}{42277273}
   \\ & &+\frac{28649412724 \sqrt{\frac{2376770}{299}} t_2^3 t_4
   t_8^2}{30747192957}-\frac{328900 \sqrt{\frac{1147}{323}} t_3^2 t_4 t_8^2}{12943}
    -\frac{1973400   \sqrt{\frac{1147}{323}} t_2 t_3 t_5 t_8^2}{12943}+\frac{95680 \sqrt{\frac{72105}{5235167}} t_5 t_6
   t_8^2}{1333}-\frac{62700755 \sqrt{13} t_3^3 t_4}{1675044}
    +\frac{13060761683 \sqrt{\frac{28823}{2530}}
   t_2^3 t_3 t_4}{41839668}
    \\ & &+\frac{13060761683 \sqrt{\frac{28823}{2530}} t_2^4
   t_5}{251038008}-\frac{62700755 \sqrt{13} t_2 t_3^2 t_5}{186116}+\frac{300099623
   \sqrt{\frac{152551}{174}} t_2^2 t_4 t_6}{575295435}-\frac{308009}{289} \sqrt{\frac{1885}{6752823}}
   t_3 t_5 t_6+t_3 t_4 t_7+t_2 t_5 t_7
     -\frac{791095425835 \sqrt{\frac{26354185}{483}} t_2 t_3^3
   t_8}{5474043792}
   \\ & &+\frac{1898304592 \sqrt{\frac{238}{1131}} t_2^2 t_4^2 t_8}{238350333}-\frac{23920
   \sqrt{\frac{72105}{5235167}} t_2 t_5^2 t_8}{1333}+\frac{614208439666441 \sqrt{\frac{232841}{114}}
   t_2^4 t_3 t_8}{993106359648}-\frac{47840 \sqrt{\frac{72105}{5235167}} t_3 t_4 t_5
   t_8}{1333}-\frac{13060761683 \sqrt{\frac{28823}{2530}} t_2^3 t_6 t_8}{125519004}+\frac{62700755
   \sqrt{13} t_3^2 t_6 t_8}{279174}
     \\ & &+t_6 t_7 t_8+\frac{68219307474560 \sqrt{\frac{7590}{50141}} t_2^2
   t_8^4 Z}{24488871563}-\frac{8610257920 \sqrt{\frac{2542}{247}} t_2 t_4 t_8^3
   Z}{869250031}+\frac{1406268931799 \sqrt{\frac{29}{13123110}} t_2^2 t_3^2 Z}{971244}+\frac{2037598
   \sqrt{\frac{7714}{49335}} t_6^2 Z}{193545}+\frac{313100288 \sqrt{\frac{336490}{1131}} t_4^2 t_8^2
   Z}{23186739199}
     \\ & &-\frac{11662696 \sqrt{\frac{62201810}{483}} t_2^2 t_3 t_8^2 Z}{1370109}+\frac{992
   \sqrt{51578} t_2 t_6 t_8^2 Z}{1677}-\frac{41478961 \sqrt{\frac{48298}{13}} t_2 t_3 t_6
   Z}{13354605}+\frac{992 \sqrt{51578} t_2 t_3 t_4 t_8 Z}{5031}-\frac{2432}{39}
   \sqrt{\frac{5474}{5488395}} t_4 t_6 t_8 Z
\end{eqnarray*}}
Where $Z$ is a solution of the quadratic  equation
\footnotesize{\begin{eqnarray*}
0&=&Z^2(6087380819159303391005603325000 \sqrt{1168090})+ Z (199260374254874324273405781024768000 \sqrt{749772699} t_1^5\\ & &-1505637451640394074610851443200000
   \sqrt{157777102} t_8 t_1^2-7430529091408703313576581198160000 \sqrt{9570} t_2 t_1)\\ & &
+326963614070713431645521199551350841540608 \sqrt{1168090}
   t_1^{10}-174838428297999356567466344499354009600 \sqrt{1604204745} t_8
   t_1^7\\ & &-1711334762076397518024060769393648435200 \sqrt{39767} t_2
   t_1^6+9960759444454760229501745087395840000 \sqrt{1168090} t_8^2 t_1^4\\ & &+9556665279777336201987345840000000 \sqrt{4109435330} t_3
   t_1^4+85101765451358857061724864000000 \sqrt{316110410} t_4
   t_1^3\\ & &+616317099570460802320908222401280000 \sqrt{1292646} t_2 t_8
   t_1^3+7240507609750674488274075585786000 \sqrt{1168090} t_2^2 t_1^2\\ & &+17177727076618842365551920000000
   \sqrt{1604204745} t_8^3 t_1+210257920682101729710627354000000 \sqrt{3059} t_5
   t_1\\ & &-147592310295208496926615407000000 \sqrt{115552965} t_3 t_8 t_1+223400586725581166350207258500000
   \sqrt{39767} t_2 t_8^2\\ & &-35432800382670021674509694653125 \sqrt{3139339} t_2
   t_3+248545233845734676184486652500 \sqrt{123400365} t_6-16850054997131291207672400000
   \sqrt{1502188545} t_4 t_8
\end{eqnarray*}}
\end{landscape}

\noindent{\bf Acknowledgments.}

 The author thanks  B. Dubrovin for useful discussions. This work was done primarily during the author  postdoctoral fellowship at the Abdus Salam International Centre for Theoretical Physics (ICTP). Italy.

\end{document}